\documentclass[10pt]{amsart}
\usepackage{amsmath,amsthm,amssymb}
\usepackage{epsfig}
\usepackage[margin=1.5in]{geometry}

\begin{document}

%%%%%%%%%%%%%%%%%%%%%%%%%%%%%%%%%%%%%%%%%%%%%%%%%%%%%%%%%%%%%%%%%%%%%%%%%%%%%%%%%%%%%
%  Group Notation
%%%%%%%%%%%%%%%%%%%%%%%%%%%%%%%%%%%%%%%%%%%%%%%%%%%%%%%%%%%%%%%%%%%%%%%%%%%%%%%%%%%%%

\newcommand{\PGL} {\Pj\Gl_2(\R)}                           %PGL_2(R)
\newcommand{\PGLC} {\Pj\Gl_2(\Cx)}                         %PGL_2(C)
\newcommand{\RP} {\R\Pj^1}                                 %RP^1
\newcommand{\CP} {\Cx\Pj^1}                                %CP^1
\newcommand{\Cx} {{\mathbb C}}                             %complex
\newcommand{\R} {{\mathbb R}}                              %reals
\newcommand{\I} {{\mathbb I}}                              %interval
\newcommand{\Z} {{\mathbb Z}}                              %integers
\newcommand{\Pj} {{\mathbb P}}                             %projective space
\newcommand{\T} {{\mathbb T}}                              %torus
\newcommand{\Sg} {{\mathbb S}}                             %symmetric group
\newcommand{\Gl} {{\rm Gl}}                                %g.linear group

\newcommand{\suchthat} {\:\: | \:\:}
\newcommand{\ore} {\ \ {\it or} \ \ }
\newcommand{\oand} {\ \ {\it and} \ \ }

%%%%%%%%%%%%%%%%%%%%%%%%%%%
% Configuration spaces
%%%%%%%%%%%%%%%%%%%%%%%%%%%

\newcommand{\oM} [1] {\ensuremath{{\mathcal M}_{0,#1}(\R)}}              %open
\newcommand{\M} [1] {\ensuremath{{\overline{\mathcal M}}{_{0, #1}(\R)}}}    %compact DMK

%%%%%%%%%%%%%%%%%%%%%%%%%%%
% Hyperplane commands
%%%%%%%%%%%%%%%%%%%%%%%%%%%

\newcommand{\C} [1]  {\mathcal C{#1}}                                      %cox complex
\newcommand{\Cp} [1] {{\Pj\C{}}({#1})}                                     %Projective cox complex
\newcommand{\Cm} [1] {\C{} ({#1})_{\#}}                                    %min blowup of cox complex
\newcommand{\Cpm} [1] {\Pj{}\Cm{#1}}                                       %min blowup of Projective cox complex
\newcommand{\Min} [1] {\textrm{Min}(\C{#1})}                               %minimal building set

\newcommand{\Hs} [1] {{\mathcal H}{#1}}                                    %hyperplanes stabilize

%%%%%%%%%%%%%%%%%%%%%%%%%%%
% Some Algebraic Geometry
%%%%%%%%%%%%%%%%%%%%%%%%%%%

\newcommand{\al}{\alpha}
\newcommand{\be}{\beta}
\newcommand{\ga}{\gamma}
\newcommand{\stp}{\perp}

\newcommand{\TS} [2] {{T_{#1}(#2)}}                                        %tangent space (bundle) of #2 at #1
\newcommand{\NS} [2] {{N_{#1}(#2)}}                                        %normal space (bundle)of #2 at #1

%%%%%%%%%%%%%%%%%%%%%%%%%%%
% Poset/Tubing commands
%%%%%%%%%%%%%%%%%%%%%%%%%%%

\newcommand{\Tubeset} {\mathfrak{T}}                     %set of all tubings
\newcommand{\Tubing} {\mathsf{T}}                        %a particular tubing
\newcommand{\tubeleq} {\prec}                            %tubes partial ordering
\newcommand{\Tleq} {\tubeleq}                            %tubings partial ordering
\newcommand{\Potset}[1] {(\mathfrak{T}_{#1},\prec_{#1})} %tube poset with order
\newcommand{\Pow} {\Omega}                               %Power set

\newcommand{\Pol} {\mathcal{P}}                          %Polytope
\newcommand{\Cox} {\Gamma}                               %coxeter graph
\newcommand{\Coxrec}[1] {\Gamma^*_{#1}}                  %reconnected complement
\newcommand{\PG}[1] {\mathcal{P}\Gamma_{#1}}             %polytope of gamma
\newcommand{\PGre}[1] {\mathcal{P}\Gamma^*_{#1}}         %polytope of gammareconnected

\newcommand{\dg}{\ \includegraphics{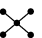} \ }                   %d5 Coxeter graph picture
\newcommand{\df}{\ \includegraphics{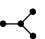} \ }                   %d4 Coxeter graph picture

%%%%%%%%%%%%%%%%%%%%%%%%%%%%%%%%%%%%%%%%%%%%%%%%%%%%%%%%%%%%%%%%%%%%%%%%%%%%%%%%%%%%%%%%%%%
%
%  paper formatting
%
%%%%%%%%%%%%%%%%%%%%%%%%%%%%%%%%%%%%%%%%%%%%%%%%%%%%%%%%%%%%%%%%%%%%%%%%%%%%%%%%%%%%%%%%%%%

\theoremstyle{plain}
\newtheorem{thm}{Theorem}[section]
\newtheorem{prop}[thm]{Proposition}
\newtheorem{cor}[thm]{Corollary}
\newtheorem{lem}[thm]{Lemma}
\newtheorem{conj}[thm]{Conjecture}

\theoremstyle{definition}
\newtheorem{defn}[thm]{Definition}
\newtheorem{exmp}[thm]{Example}

\theoremstyle{remark}
\newtheorem*{rem}{Remark}
\newtheorem*{hnote}{Historical Note}
\newtheorem*{nota}{Notation}
\newtheorem*{ack}{Acknowledgments}
\numberwithin{equation}{section}

\title {Coxeter Complexes and Graph-Associahedra}

\thanks{Authors were partially supported by NSF grants DMS-9820570 and CARGO DMS-0310354.}

\author[M.\ Carr]{Michael Carr}
\address{M.\ Carr: University of Michigan, Ann Arbor, MI 48109}
\email{mpcarr@umich.edu}

\author[S.\ Devadoss]{Satyan L.\ Devadoss}
\address{S.\ Devadoss: Williams College, Williamstown, MA 01267}
\email{satyan.devadoss@williams.edu}

\begin{abstract}
Given a graph $\Cox$, we construct a simple, convex polytope, dubbed \emph{graph-associahedra}, whose face poset is based on the connected subgraphs of $\Cox$.  This provides a natural generalization of the Stasheff associahedron and the Bott-Taubes cyclohedron.  Moreover, we show that for any simplicial Coxeter system, the minimal blow-ups of its associated Coxeter complex has a tiling by graph-associahedra.  The geometric and combinatorial properties of the complex as well as of the polyhedra are given.   These spaces are natural generalizations of the Deligne-Knudsen-Mumford compactification of the real moduli space of curves.
\end {abstract}


\maketitle

%%%%%%%%%%%%%%%%%%%%%%%%%%%%%%%%%%%%%%%%%%%%%%%%%%%%%%%%%%%%%%%%%%%%%%%%%%%%%%%%%%%%%
%NOTE: The \baselineskip=15pt below is the minimum needed to make the paper look decent.
%      There are many superscripts {\tilde, \overline} which cause erratic
%      spaces between sentences. \baselineskip corrects these problems.  Also see the
%      \baselineskip=12pt by References.
%%%%%%%%%%%%%%%%%%%%%%%%%%%%%%%%%%%%%%%%%%%%%%%%%%%%%%%%%%%%%%%%%%%%%%%%%%%%%%%%%%%%%

\baselineskip=15pt

%%%%%%%%%%%%%%%%%%%%%%%%%%%%%%%%%%%%%%%%%%%%%%%%%%%%%%%%%%%%%%%%%%%%%%%%%%%%%%%%%%%%%
%
%                Introduction
%
%%%%%%%%%%%%%%%%%%%%%%%%%%%%%%%%%%%%%%%%%%%%%%%%%%%%%%%%%%%%%%%%%%%%%%%%%%%%%%%%%%%%%
\section{Introduction}

The Deligne-Knudsen-Mumford compactification of the real moduli space of curves \break \M{n} appears in many areas, from operads \cite{dev, sta2} , to combinatorics \cite{dev3, kap1}, to group theory \cite{djs, djs2}.  One reason for this is an intrinsic tiling of \M{n} by the \emph{associahedron}, the Stasheff polytope \cite{sta}.  The motivation for this work comes from a remarkable fact, first noticed by Kapranov, involving Coxeter complexes:   Blowing up certain faces of the Coxeter complex of type $A$ yields a double cover of \M{n}.  Extending this to the Coxeter complex of affine type $\widetilde A$ results in a moduli space tessellated by the \emph{cyclohedron} \cite{dev2}, the Bott-Taubes polytope associated to knot invariants \cite{bt}.  Davis et al.\  have shown these spaces to be aspherical, where all the homotopy properties are completely encapsulated in their fundamental groups \cite{djs}.
This paper looks at analogues of \M{n} for all simplicial Coxeter groups $W$, which we denote as $\Cm{W}$.

Section~\ref{s:gassoc} begins with the study of graph-associahedra.  For any graph $\Cox$, we construct a simple, convex polytope whose face poset is based on the connected subgraphs of $\Cox$ (Theorem~\ref{t:faceposet}).  This provides a natural generalization of the associahedron and the  cyclohedron.  Some combinatorial properties of this polytope are also explored (Theorem~\ref{t:codim1}).

Section~\ref{s:tiling} provides the background of Coxeter complexes and proves that graph-associa-hedra tile $\Cm{W}$ (Theorem~\ref{t:tile}).  A gluing map of these polytopes is also provided (Theorem~\ref{t:glue}).  Section~\ref{s:geometry} finishes by looking at the geometry of $\Cm{W}$.  In particular, we show that each blown-up cell of $\Cm{W}$ resolves into a product of lower-dimensional blown-up Coxeter complexes (Theorem~\ref{t:blow_ups}).

%%%%%%%%%%%%%%%%%%%%%%%%%%%%%%%%%%%%%%%%%%%%%%%%%%%%%%%%%%%%%%%%%%%%%%%%%%%%%%%%%%%%%
%
%                Constructing Graph-Associahedra
%
%%%%%%%%%%%%%%%%%%%%%%%%%%%%%%%%%%%%%%%%%%%%%%%%%%%%%%%%%%%%%%%%%%%%%%%%%%%%%%%%%%%%%
\section{Constructing Graph-Associahedra}
\label{s:gassoc}
\subsection{}

The motivating example will be the associahedron.

\begin{defn}
Let ${\mathfrak A}(n)$ be the poset of bracketings of a path with $n$ nodes, ordered such that $a \prec a'$ if $a$ is obtained from $a'$ by adding new brackets.  The \emph{associahedron} $K_n$ is a convex polytope of dimension $n-2$ whose face poset is isomorphic to ${\mathfrak A}(n)$.
\end{defn}

\begin{figure}[h]
\includegraphics {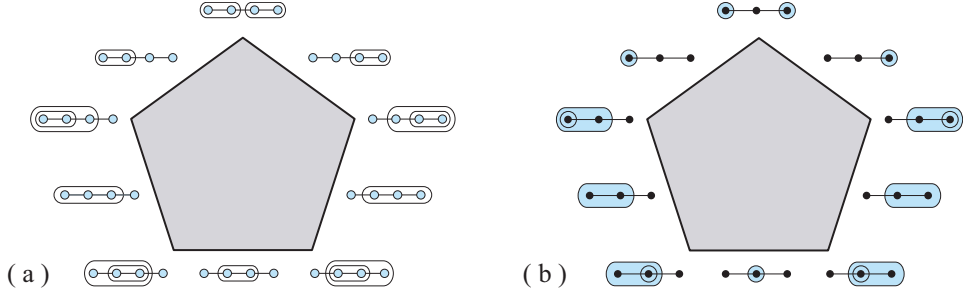}
\caption{Associahedron $K_4$ labeled with (a) bracketings and (b) tubings.}
\label{f:k4pt}
\end{figure}

The associahedron $K_n$ was originally defined by Stasheff for use in homotopy theory in connection with associativity properties of $H$-spaces \cite[Section 2]{sta}.  The construction of the polytope $K_n$ is given by Lee \cite{lee} and Haiman (unpublished). The vertices of $K_n$ are enumerated by the Catalan numbers.  Figure~\ref{f:k4pt}(a) shows the $2$-dimensional $K_4$ as the pentagon. Each edge of $K_4$ has one set of brackets, whereas each vertex has two. Note that Figure~\ref{f:coxblow}(a) shows $\Cm{A_3}$ tiled by 24 \,$K_4$ pentagons.  We give an alternate definition of $K_n$ with respect to \emph{tubings}.

\begin{defn}
Let $\Cox$ be a graph.  A \emph{tube} is a proper nonempty set of nodes of $\Cox$ whose induced graph is a proper, connected subgraph of $\Cox$.  There are three ways that two tubes $t_1$ and $t_2$ may interact on the graph.
\begin{enumerate}
  \item Tubes are \emph{nested} if  $t_1 \subset t_2$.
  \item Tubes \emph{intersect} if $t_1 \cap t_2 \neq \emptyset$ and $t_1 \not\subset t_2$ and $t_2 \not\subset t_1$.
  \item Tubes are \emph{adjacent} if $t_1 \cap t_2 = \emptyset$ and $t_1 \bigcup t_2$ is a tube in $\Cox$.
\end{enumerate}
Tubes are \emph{compatible} if they do not intersect and they are not adjacent.  A \emph{tubing} $T$ of $\Cox$ is a set of tubes of $\Cox$ such that every pair of tubes in $T$ is compatible.  A \emph{$k$-tubing} is a tubing with $k$ tubes.
\end{defn}

\begin{lem}
Let $\Cox$ be a path with $n-1$ nodes.  The face poset of $K_n$ is isomorphic to the poset of all valid tubings of \ $\Cox$, ordered such that tubings \ $T \prec T'$  if $T$ is obtained from $T'$ by adding tubes.
\end{lem}

\noindent Figure~\ref{f:k4pt}(b) shows the faces of associahedron $K_4$ labeled with tubings.  The proof of the lemma is based on a trivial bijection between bracketings and tubings on paths.

%%%%%%%%%%%%%%%%%%%%%%%%%%%%%%%%%%%%%%%%%%%%%%%%%%%%%%%%%%%%%%%%%%%%%%%%%%%%%%%%%%%%%

\subsection{}
For a graph $\Cox$ with $n$ nodes, let $\triangle_{\Cox}$ be the $n-1$ simplex in which each facet (codimension $1$ face) corresponds to a particular node.  Each proper subset of nodes of $\Cox$ corresponds to a unique face of $\triangle_{\Cox}$, defined by the intersection of the faces associated to those nodes. The empty set corresponds to the face which is the entire polytope $\triangle_{\Cox}$.

\begin{defn}
For a given graph $\Cox$, truncate faces of $\triangle _{\Cox}$ which correspond to $1$-tubings in increasing order of dimension. The resulting polytope $\PG{}$ is the \emph{graph-associahedron}.
\end{defn}

\noindent This definition is well-defined:  Theorem~\ref{t:faceposet} below guarantees that truncating any ordering of faces of the same dimension produces the same poset/polytope.  Note also that $\PG{}$ is a simple, convex polytope.

\begin{figure}[h]
\includegraphics {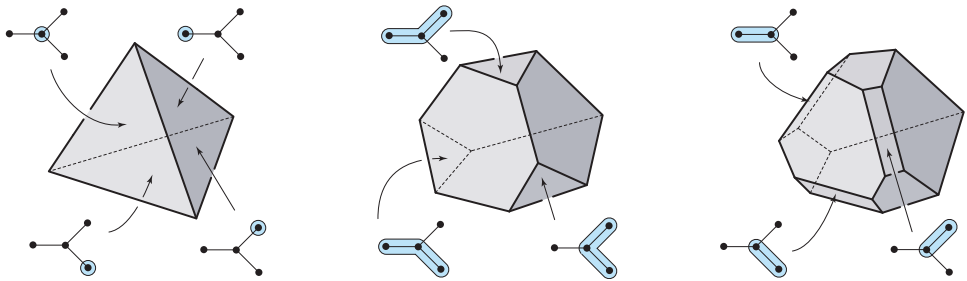}
\caption{Iterated truncations of the $3$-simplex based on an underlying graph.}
\label{f:d4}
\end{figure}

\begin{exmp}
Figure \ref{f:d4} shows a $3$-simplex tetrahedron truncated according to a graph.  The facets of $\Pol ( \df )$ are labeled with $1$-tubings.  One can verify that the edges correspond to all possible 2-tubings and the vertices to $3$-tubings.
\end{exmp}

\begin{thm}
$\PG{}$ is a simple, convex polytope whose face poset is isomorphic to set of valid tubings of \,$\Cox$, ordered such that $T \prec T'$ if $T$ is obtained from $T'$ by adding tubes.
\label{t:faceposet}
\end{thm}

\noindent The proof of this theorem is given at the end of the section.  Note that simplicity and convexity of $\PG{}$ follows from its construction.  Stasheff and Schnider \cite[Appendix B]{sta2} proved the following motivating examples.  They follow immediately from Theorem~\ref{t:faceposet}.

\begin{cor}
When $\Cox$ is a path with $n-1$ nodes, $\PG{}$ is the associahedron $K_n$.  When $\Cox$ is a cycle with $n-1$ nodes, $\PG{}$ is the cyclohedron $W_n$.
\end{cor}

%%%%%%%%%%%%%%%%%%%%%%%%%%%%%%%%%%%%%%%%%%%%%%%%%%%%%%%%%%%%%%%%%%%%%%%%%%%%%%%%%%%%%

\subsection{}
For a given tube $t$ and a graph $\Cox$, let $\Cox_{t}$ denote the induced subgraph on the graph $\Cox$.  By abuse of notation, we sometimes refer to $\Cox_{t}$ as a tube.

\begin{defn}
Given a graph $\Cox$ and a tube $t$, construct a new graph $\Coxrec{t}$ called the \emph{reconnected complement}: If $V$ is the set of nodes of $\Cox$, then $V-t$ is the set of nodes of $\Coxrec{t}$.  There is an edge between nodes $a$ and $b$ in $\Coxrec{t}$ if either $\{a,b\}$ or $\{a,b\} \cup t$ is connected in $\Cox$.
\end{defn}

\begin{figure}[h]
\includegraphics {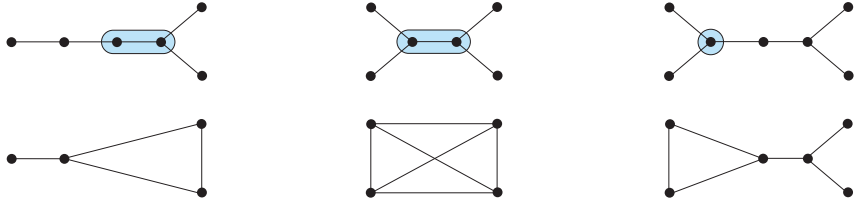}
\caption{Examples of $1$-tubings and their reconnected complements.} \label{f:recon}
\end{figure}

\noindent Figure~\ref{f:recon} illustrates some examples of $1$-tubings on graphs along with their reconnected complements.

\begin{thm}
The facets of \ $\PG{}$ correspond to the set of \ $1$-tubings on $\Cox$.  In particular, the facet associated to a $1$-tubing $\{t\}$ is combinatorially equivalent to $\PG{t} \times \PGre{t}$.
\label{t:codim1}
\end{thm}

\begin{proof}
We know from Theorem~\ref{t:faceposet} that a facet of $\PG{}$ is given by a $1$-tubing $\{t\}$.  The faces contained in this facet are the tubings $T$ of $\Cox$ that contain $t$.  Now if $t_i \subset t$ is a tube of $\Cox_t$ then it is also a tube of $\Cox$.  Consider the map
$$\rho : \{\text{tubes of $\Coxrec{t}$}\} \to \{\text{tubes of $\Cox$ containing $t$}\}$$
where
$$\rho (t') =
\begin{cases}
\ \ t' \cup t & \hspace{.5in} \text{if $t' \cup t$ is connected in $\Cox$} \\
\ \ t' & \hspace{.5in} \text{otherwise}.
\end{cases}$$
Note that $\rho$ is a bijection and it preserves the validity of tubings.  That is, two tubes $t_1$ and $t_2$ are compatible in $\Coxrec{t}$ if and only if $\rho(t_1)$ and $\rho(t_2)$ are compatible.  Define the natural map
$$\widehat \rho : \{ \text{tubings on $\Coxrec{t}$} \} \times \{ \text{tubings on $\Cox_{t}$} \} \to
\{ \text{tubings on $\Cox$} \}$$ where
$$\widehat \rho(T_i \times T_j) \ = \ \{t\} \ \cup \bigcup_{t_i \in T_i} \{\rho (t_i)\} \ \cup \bigcup_{t_j\in T_j} \{t_j\}.$$
It is straightforward to show that this is an isomorphism of posets.
\end{proof}

\begin{figure}[h]
%\resizebox{\textwidth}{!}{
\includegraphics{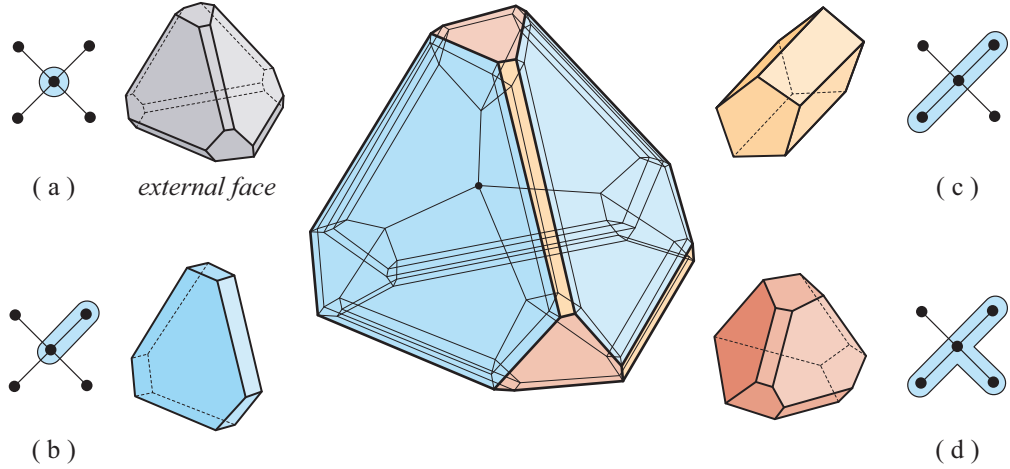}
\caption{The Schlegel diagrams of a $4$-polytope along with its four types of facets.}
\label{f:d5}
\end{figure}

\begin{exmp}
Figure~\ref{f:d5} shows the Schlegel diagram of the $4$-dimensional polytope $\Pol ( \dg )$.  It is obtained from the $4$-simplex by first truncating four vertices, each of which become a $3$-dimensional facet, as depicted in Figure~\ref{f:d5}(d) along with its $1$-tubing.  Then six edges are truncated, becoming facets of type Figure~\ref{f:d5}(c); note that Theorem~\ref{t:codim1} shows the structure of the facet to be the product of the associahedron $K_4$ of Figure~\ref{f:k4pt}(b) and an interval.  Finally four $2$-faces of the original $4$-simplex are truncated to result in the polytope of Figure~\ref{f:d5}(b); this is the product of the cyclohedron $W_3$ (hexagon) and an interval.  Four of the original five facets of the $4$-simplex have become the polyhedron of Figure~\ref{f:d5}(d), whereas the fifth (external) facet is the $3$-dimensional \emph{permutohedron}, as shown in Figure~\ref{f:d5}(a).
\end{exmp}

%%%%%%%%%%%%%%%%%%%%%%%%%%%%%%%%%%%%%%%%%%%%%%%%%%%%%%%%%%%%%%%%%%%%%%%%%%%%%%%%%%%%%

\subsection{}
The remaining section is devoted to the proof of Theorem~\ref{t:faceposet}, which follows directly from Lemmas~\ref{l:kequalsf} and \ref{l:truncate} below.  First we must define a poset operation analogous to truncation.

We define an initial partial ordering $\tubeleq_0$ on tubes by saying that $t_{i} \tubeleq_0 t_{j}$ if and only if $t_{i} \subset t_{j}$.  We also define a partial ordering on a set of tubings $\Tubeset$ induced by any partial ordering of tubes of $\Cox$: Given tubings $T_I, \, T_J \in \Tubeset$, then $T_I \Tleq T_J$ if and only if for all $t_j \in T_J$, there exists $t_i$ such that $t_j \tubeleq t_i \in T_I$.
%Abusing notation, we understand by context that $\Tleq$ refers to both ordering of tubes
% and the induced ordering on tubings.
We write this partially ordered set of tubings as $\Potset{}$. Note that $\triangle_{\Cox}$ is isomorphic to $\Potset{0}$: the set of nonnested tubings of $\Cox$ with order induced by $\tubeleq_0$.
% We will write $\Tleq_i$  to indicate the ordering induced by $\tubeleq_i$.

\begin{defn} Given a poset of tubings $(\Tubeset,\tubeleq)$, we can produce a set $(\Tubeset', \tubeleq ')$ by \emph{promoting} the tube $t_{*}$.  Let
$$\Tubeset_* = \{ T\cup \{t_{*}\} \suchthat T \in \Tubeset \ \text{and} \ T\cup \{t_{*}\} \ {\text{is a valid tubing of}}\ \Cox \}$$
and let $\Tubeset'=\Tubeset\cup \Tubeset_*$.  Let $\tubeleq'$ be defined so that $t_*$ is incomparable to any other tube, and for any tubes $t_{a}, t_{b}$ not equal to $t_{*}$, let $t_{a} \tubeleq' t_{b}$ if and only if $t_{a} \tubeleq t_{b}$.
\end{defn}

Let $\{t_{i}\}$ be the set of tubes in $\Cox$ for $ i \in \{ 1, \dotsc, k\}$ ordered in decreasing size.  Notice these correspond to the faces of $\triangle_{\Cox}$ in increasing order of dimension. Let $\Potset{i}$ be the resulting set after consecutively promoting the tubes $t_{1}, \dotsc, t_{i}$ in $\Potset{0}$.  The following two lemmas explicitly define the tubings and the ordering of $\Potset{i}$.  Both are trivial inductions from the definition of promotion.

%%%%%%%%%%%%%%%%%%%%%%  LEMMA    %%%%%%%%%%%%%%%%%%

\begin{lem}
$\Tubeset_i$ is the set of all valid tubings of the form $T_0 \bigcup_{j=1}^m \{t_{q_j}\}$ where $T_0 \in \Tubeset_0$ and $\{q_j\}\subseteq\{1,..., i\}$. \label{l:genform}
\end{lem}

%%%%%%%%%%%%%%%%%%%%%%  LEMMA    %%%%%%%%%%%%%%%%%%

\begin{lem}
If $a$ or $b$ is less than or equal to $i$, then $t_a\tubeleq_{i}t_b$ if and only if $a=b$.  If both $a$ and $b$ are greater than $i$, then $t_a\tubeleq_{i}t_b$ if and only if $t_a\tubeleq_{0}t_b$. \label{l:genorder}
\end{lem}

\noindent As a special case we can state the following:

%%%%%%%%%%%%%%%%%%%%%%  LEMMA    %%%%%%%%%%%%%%%%%%

\begin{lem}
$\Potset{k}$ is isomorphic as a poset to the set of tubings of \, $\Cox$, ordered such that $T\prec T'$ if and only if $T$ can be obtained by adding tubes to $T'$. \label{l:kequalsf}
\end{lem}

\begin{proof}
Applying Lemma~\ref{l:genform} to the case $i=k$ shows $\Tubeset_k$ is the set of all tubings of $\Cox$.  Lemma~\ref{l:genorder} shows that $T \tubeleq_k T'$ if and only if $T \supset T'$.
\end{proof}

\begin{figure}[h]
\resizebox{\textwidth}{!}{\includegraphics {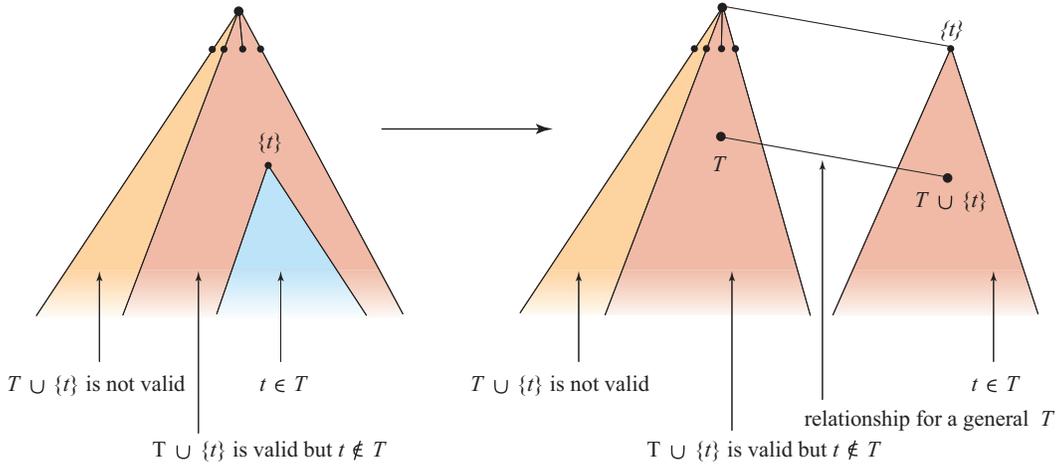}}
\caption{A sketch of the poset lattice before and after promotion of tube $\{t\}$.  Regions shaded with like colors are isomorphic as posets.} \label{f:promotion}
\end{figure}

\noindent The only step that remains is to show the equivalence of promotion to truncation when performed in this order.
The following lemma accomplishes this.
%The following lemmas accomplish this.

%%%%%%%%%%%%%%%%%%%%%%  LEMMA    %%%%%%%%%%%%%%%%%%

\begin{lem}
Let $f_{i}$ be a face of $\triangle _{\Cox}$ corresponding to the tube $t_i$.  Let $\Pol_{i}$ be the polytope created by consecutively truncating faces $f_{1},...,f_{i}$ of $\triangle _{\Cox}$. Then $\Potset{k}\cong P_{k}$. \label{l:truncate}
\end{lem}

\begin{proof}
For consistency, we refer to $\triangle _{\Cox}$ by $\Pol_0$. Since $\Pol_0$ is convex, so is $\Pol_{i}$.  Thus we may define these polytopes as intersections of halfspaces. Denote the hyperplane that defines the halfspace $H_a^+$ by $H_a$.  If $X$ is the halfspace set for a polytope $\Pol$ then there is a natural poset map
$$\Psi : \Pol \to \Pow(X)^{op} : f \mapsto X_f$$
where $\Pow(X)^{op}$ is the set of subsets of $X$ ordered under reverse inclusion and $X_f$ is the subset such that $f=\Pol\cap \bigcap_{a\in X_f}H_a$. Note that $\Psi$ is an injection with its image as all the sets $X'$ such that $\Pol\cap \bigcap_{a\in X'}H_a$ in nonempty. By truncating $\Pol$ at $f_*$, a new halfspace $H_*^+$ is added with the following properties:
\begin{enumerate}
\item A vertex of $\Pol$ is in $H_*^+$ if and only if it is not in $f_*$. \item No vertices of $\Pol$ are in $H_*$.
\end{enumerate}
This produces the truncated polytope $\Pol_*=H_*^+ \ \cap \ \bigcap_{a\in X}H_a^+$.  Let  $\Pol_0$ be defined by $\bigcap_{a\in X_0}H_a^+$  where $X_0$ is the set of indices for the defining halfspaces.   Let $H_i^+$ be the halfspace with which we intersect $\Pol_{i-1}$ to truncate $f_i$.   The halfspace set for $\Pol_i$ is
$$X_i \ = \ X_{i-1} \cup \{i\} \ = \ X_0 \cup \{1, \dotsc, i\}.$$
We define the map $\Psi_i :\Pol_i \to \Pow(X_i)^{op}$ which takes a face of $\Pol_i$ to the set of hyperplanes that contain it.

We now produce an order preserving injection $\Phi_i$ from $\Tubeset_i$ to $\Pow(X_i)^{op}$.  Let $\phi_0$ be the map from tubes of $\Cox$ to $\Pow(X_0)$ that takes a tube $t_i$ to $\Psi_0(f_i)$.  Define
$$\phi_i (t_j) =
\begin{cases}
\ \ \{j\} & \ \text{if $j \leq i$}, \\
\ \phi_0(t_j) & \ \text{if $j >  i$}.
\end{cases}$$
This allows us to define a new map
$$\Phi_i : \Tubeset_i \to \Pow(X_i)^{op} : T_J \mapsto \bigcup_{t_j \in T_J} \phi_i(t_j).$$

\noindent It follows from the definition that this is an order preserving injection.
An induction argument %Lemma~\ref{l:lastcase} below
shows that $\Phi_i(\Tubeset_i) = \Psi_i (\Pol_i)$.
Since $\Psi_i$ and $ \Phi_i$ are order preserving and injective, we have that $\Psi_i^{-1} \circ \Phi_i:\Tubeset_i \to \Pol_i$ is an isomorphism of posets.
\end{proof}

\section{Tiling Coxeter Complexes}
\label{s:tiling}
\subsection{}

We begin with some standard facts and definitions about Coxeter systems.  Most of the background used here can be found in Bourbaki \cite{bou} and Brown \cite{bro}.
\begin{defn}
Given a finite set $S$, a \emph{Coxeter group} $W$ is given by the presentation
$$W \ = \ \langle \ s_i \in S \ \ | \ \ s_i^2 = 1, \ (s_i s_j)^{m_{ij}} = 1 \ \rangle \, ,$$
where $m_{ij} = m_{ji}$ and $2 \leq m_{ij} \leq \infty$.
\end{defn}
\noindent Associated to any Coxeter system $(W, S)$ is its \emph{Coxeter graph} $\Cox_W$:  Each node represents an element of $S$, where two nodes $s_i, s_j$ determine an edge if and only if $m_{ij} \geq 3$.  A Coxeter group is \emph{irreducible} if its Coxeter graph is connected and it is \emph{locally finite} if either $W$ is finite or each proper subset of $S$ generates a finite group.  A Coxeter group is \emph{simplicial} if it is irreducible and locally finite. The classification of simplicial Coxeter groups and their Coxeter graphs are well-known \cite[Chapter 6]{bou}.  Unless stated otherwise, the Coxeter groups discussed below are assumed to be simplicial.

Every simplicial Coxeter group has a realization as a group generated by reflections acting faithfully on a variety \cite[Chapter 3]{bro}.  The geometry of the variety is either spherical, Euclidean, or hyperbolic, depending on the group.  Every conjugate of a generator $s_i$ acts on the variety as a reflection in some hyperplane, dividing the variety into simplicial chambers.  This variety, along with its cellulation is the \emph{Coxeter complex} corresponding to $W$, denoted $\C{W}$.  The hyperplanes associated to the generators $s_i$ of $W$ all border a single chamber, called the \emph{fundamental chamber} of $\C{W}$.  The $W$-action on the chambers of $\C{W}$ is transitive, and thus we may associate an element of $W$ to each chamber; generally, the identity is associated to the fundamental chamber.

%\begin{rem}
%A Coxeter complex of a spherical or Euclidean type is an algebraic variety in an obvious way.  The equations that define hyperplanes are linear and affine equations, respectively.  In the case of a hyperbolic complex, the hyperbolic space (minus the origin) maps conformally to the set $x_1>0$ in $\R^n$.  Thus taking $\R^n$ and subtracting the $x_1=0$ plane gives two identical copies of hyperbolic space.  This is in fact a variety, and the hyperplanes are defined by the equations of spheres centered on $x_1=0$ or planes orthogonal to $x_1=0$.  We will abuse notation and use the term \emph{Coxeter complex} to refer to the associated variety.  The results that follow can easily be applied to the actual Coxeter complex.
%\end{rem}

\begin{nota}
For a spherical Coxeter complex $\C{W}$, we define the \emph{projective Coxeter complex} $\Cp{W}$ to be $\C{W}$ with antipodal points on the sphere identified.  These complexes arise naturally in blow-ups, as shown in Theorem~\ref{t:blow_ups}.
\end{nota}

\begin{figure}[h]
\includegraphics{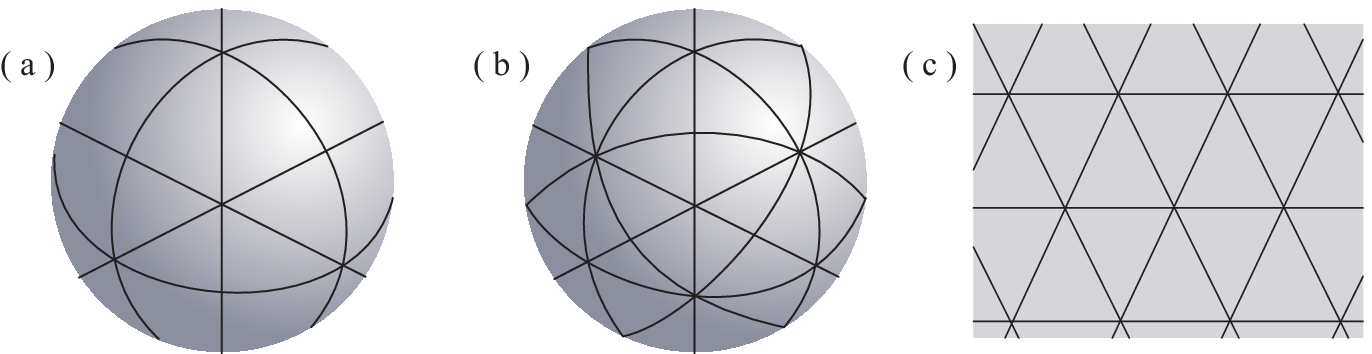}
\caption{Coxeter complexes $\C{A_3}$, $\C{B_3}$, and $\C{{\widetilde A_2}}$.}
\label{f:cox}
\end{figure}

\begin{exmp}
The Coxeter group of type $A_n$ has $n$ generators, and $m_{ij}=3$ if $i=j\pm 1$ and $2$ otherwise.  Thus $A_n$ is isomorphic to the symmetric group $\Sg_{n+1}$ and acts on the intersection of the unit sphere in $\R^{n+1}$ with the hyperplane $x_1+x_2+\cdots +x_{n+1}=0$.  Each $s_i$ is the reflection in the plane $x_i=x_{i+1}$.  Figure~\ref{f:cox}(a) shows the Coxeter complex $\C{A_3}$, a $2$-sphere cut into $24$ triangles.

The $B_n$ Coxeter group has $n$ generators with the same $m_{ij}$ as $A_n$ except that $m_{12}=4$.  The group $B_n$ is the symmetry group of the $n$-cube, and acts on the unit sphere in $\R^n$.  Each generator $s_i$ is a reflection in the hyperplane $x_{i-1}=x_i$, except $s_1$ which is the reflection in $x_1=0$.  Figure~\ref{f:cox}(b) shows the Coxeter complex $\C{B_3}$, the $2$-sphere tiled by simplices.

The $\widetilde A_{n}$ Coxeter group has $n+1$ generators, with $m_{ij}=3$ if $i=j\pm 1$, and $m_{(1)(n+1)}=3$.  Every other $m_{ij}$ equals two. The group ${\widetilde A_{n}}$ acts on the hyperplane defined by $x_1+x_2+\cdots +x_{n+1}=0$ in $\R^{n+1}$. Each $s_i$ is the reflection in $x_i=x_{i+1}$, except $s_{n+1}$ which is the reflection in $x_{n+1}=x_1+1$.  Figure~\ref{f:cox}(c) shows the Coxeter complex $\C{{\widetilde A_2}}$, the plane with the corresponding hyperplanes.
\end{exmp}

%%%%%%%%%%%%%%%%%%%%%%%%%%%%%%%%%%%%%%%%%%%%%%%%%%%%%%%%%%%%%%%%%%%%%%%%%%%%%%%%%%%%%

\subsection{}

The collection of hyperplanes $\{x_i = 0 \ | \ i = 1, \ldots, n\}$ of $\R^n$ generates the \emph{coordinate} arrangement.  A crossing of hyperplanes is \emph{normal}\, if it is locally isomorphic to a coordinate arrangement.   A construction which transforms any crossing into a normal crossing involves the algebro-geometric concept of a blow-up; see Section~\ref{ss:blowup} for a definition.

A general collection of blow-ups is usually noncommutative in nature; in other words, the order in which spaces are blown up is important.  For a given arrangement, De Concini and Procesi \cite[Section 3]{dp} establish the existence (and uniqueness) of a \emph{minimal building set}, a collection of subspaces for which blow-ups commute for a given dimension, and for which every crossing in the resulting space is normal.  We denote the minimal building set of an arrangement $\mathcal A$ by $\textrm{Min}({\mathcal A})$.  Let $\al$ be an intersection of hyperplanes in an arrangement $\mathcal A$.  Denote $\Hs{\al}$ to be the set of all hyperplanes that contain $\al$.  We say $\Hs \al$ is \emph{reducible} if it is a disjoint union $\Hs \be \sqcup \Hs \ga$, where $\al = \be \cap \ga$ for intersections of hyperplanes $\be$ and $\ga$.

\begin{lem} \cite[Section 2]{dp}
\label{l:minirred}
$\al \in \textrm{Min}({\mathcal A})$ if and only if \ $\Hs \al$ is irreducible.
\end{lem}

If reflections in $\Hs \al$ generate a Coxeter group (finite reflection group), it is called the \emph{stabilizer} of $\al$ and denoted $W_\al$.
For a Coxeter complex $\C{W}$, we denote its minimal building set by $\Min{W}$.
%Note that in a Coxeter complex, $W_\al$ exists for all intersections of hyperplanes.
The relationship between the set $\Min{W}$ and the group $W$ is given by the following.

\begin{lem} \cite[Section 3]{djs}
\label{l:irreducible}
$\al \in \Min{W}$ if and only if \ $W_\al$ is irreducible.
\end{lem}

\begin{defn}
The \emph{minimal blow-up} of $\C{W}$, denoted as $\Cm{W}$, is obtained by blowing up along elements of $\Min{W}$ in \emph{increasing} order of dimension.
\end{defn}

\noindent The construct $\Cm{W}$ is well-defined:  Lemma~\ref{l:reorder} below guarantees that blowing-up any ordering of subspaces in $\Min{W}$ of the same dimension produces the same cellulation.

\begin{figure}[h]
\includegraphics {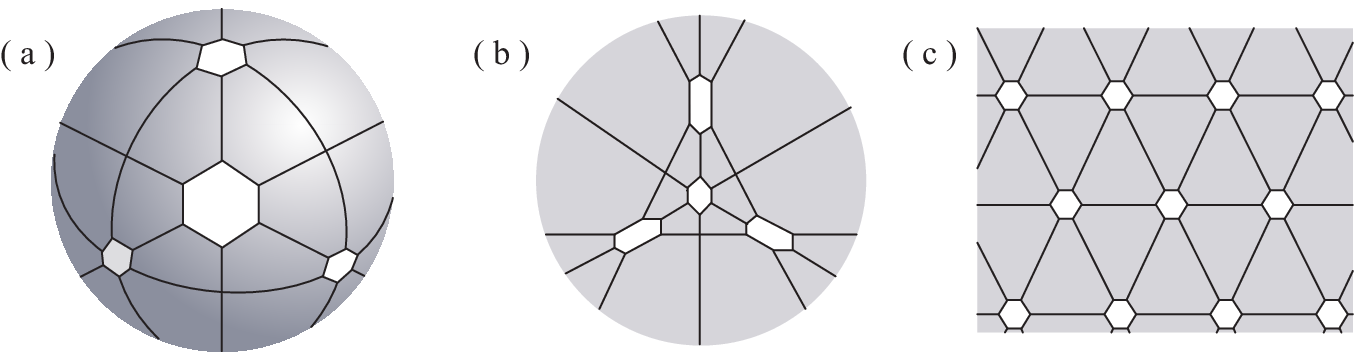}
\caption{Minimal blow-ups of (a) $\C{A_3}$, (b) $\Cp{A_3}$ and (c) $\C{\widetilde A_2}$.}
\label{f:coxblow}
\end{figure}

\begin{exmp}
Figure \ref{f:coxblow}(a) shows the blow-ups of the sphere $\C{A_3}$ of Figure~\ref{f:cox}(a) at nonnormal crossings.  Each blown up point has become a hexagon with antipodal identification and the resulting manifold is $\Cm{A_3}$.  Figure~\ref{f:gluing} shows the local structure at a blow-up, where each crossing is now normal.  The minimal blow-up of the projective Coxeter complex of type $A_3$ is shown in Figure \ref{f:coxblow}(b), with the four points blown up in $\R \Pj^2$.  Figure~\ref{f:coxblow}(c) shows the minimal blow-up of $\C{\widetilde A_2}$ of Figure~\ref{f:cox}(c).
\end{exmp}

%%%%%%%%%%%%%%%%%%%%%%%%%%%%%%%%%%%%%%%%%%%%%%%%%%%%%%%%%%%%%%%%%%%%%%%%%%%%%%%%%%%%%

\subsection{}
Given the construction of graph-associahedra above, we turn to applying them to the chambers tiling $\Cm{W}$.

\begin{thm}
Let $W$ be a simplicial Coxeter group and $\Cox_W$ be its associated Coxeter graph.  Then $\PG{W}$ is the fundamental domain for $\Cm{W}$.
\label{t:tile}
\end{thm}

\begin{proof}
It is a classic result of geometric group theory that each chamber of a simplicial Coxeter complex $\C{W}$ is a simplex.  The representation of $W$ can be chosen such that the generators correspond to the reflections through the supporting hyperplanes of a fixed chamber.  In other words, a fundamental chamber of $\C{W}$ is the simplex $\triangle_{\Cox_W}$ such that each facet of $\triangle_{\Cox_W}$ is associated to a node of $\Cox_W$.

Let $f$ be a face of  $\triangle_{\Cox_W}$ and let $\al$ be the \emph{support} of $f$, the smallest intersection of hyperplanes of $\C{W}$ containing $f$. As in the previous section, the face $f$ corresponds to a subset $S$ of the nodes of $\Cox_W$.  The nodes in $S$ represent the generators of $W$ that stabilize $\al$. These elements generate $W_\al$, and the subgraph induced by $S$ is the Coxeter graph of $W_\al$.

By Lemma~\ref{l:irreducible}, $\al$ is an element of $\Min{W}$ if and only if $W_\al$ is irreducible.  But $W_\al$ is irreducible if and only if $\Cox_{W_\al}$ is connected, that is, when the set of nodes of $\Cox_{W_\al}$ is a tube of $\Cox_W$. Note that blowing up $\al$ in $\C{W}$ truncates the face $f$ of $\triangle_{\Cox_W}$.  Thus performing minimal blow ups of $\C{W}$ is equivalent to truncating the faces of $\triangle_{\Cox_W}$ that correspond to tubes of $\Cox_W$.  By definition, the resulting polytope is $\PG{W}$.
\end{proof}

\begin{rem}
The \emph{maximal} building set is the collection of \emph{all} crossings, not just the nonnormal ones.
The fundamental chambers of the maximal blow-up of $\C{W}$ will be tiled by permutohedra, obtained by iterated truncations of all faces of the simplex.
\end{rem}

\begin{rem}
The generalized associahedra of Fomin and Zelevinsky \cite{fz} are fundamentally different than graph-associahedra.  Although both are motivated from type $A_n$ (the classical associahedra of Stasheff), they are distinct in all other cases.  For example, the cyclohedron is the generalized associahedron of type $B_n$, whereas it is the type $\widetilde A_n$ graph-associahedron.
\end{rem}

%%%%%%%%%%%%%%%%%%%%%%%%%%%%%%%%%%%%%%%%%%%%%%%%%%%%%%%%%%%%%%%%%%%%%%%%%%%%%%%%%%%%%

\subsection{}
The construction of the Coxeter complex $\C{W}$ implies a natural $W$-action.  This action, restricted to the chambers is faithful and transitive, so we can identify each chamber with the group element that takes the fundamental chamber to it. The faces of the chambers of $\C{W}$ have different types (according to their associated tubings in $\Cox_W$).  A transformation is \emph{type preserving} if it takes each face to a face of the same type.   We call the $W$-action type preserving because each $w$ induces a type preserving transformation of $\C{W}$.

We may use this action to define a $W$-action on $\Cm{W}$.  There is a hyperplane-preserving isomorphism between $\C{W}-\bigcup\Min{W}$ and $\Cm{W}-\bigcup\Min{W}$.  We define the $W$-action on $\Cm{W}$ to agree with the $W$-action on $\C{W}$ in $\Cm{W}-\bigcup\Min{W}$. We define the action on the remainder of $\Cm{W}$ by requiring that for all subvarieties $V$ of $\Cm{W}-\bigcup\Min{W}$, the action of $w$ takes the closure of $V$ to the closure of $wV$.  The $W$-action defined this way is type preserving, and the stabilizer of each hyperplane $\al$ is the group $W_\al$.

Given $\Cm{W}$, we may associate an element $s_f \in W$ to each facet $f$ of the fundamental chamber. We call this element the \emph{reflection} in that facet. If $\al$ is the hyperplane of $\Cm{W}$ that contains $f$, then $s_f$ is a reflection in $\al$. This corresponds to the reflection across $\al$ in $\C{W}$, which is the longest word in $W_{\al}$ \cite[Section 3]{bro}.  For a face $f$ of the fundamental chamber, define $W_f=\langle s_{f_i} \rangle$ and $s_f=\prod s_{f_i}$, where the $f_i$'s are the facets of the fundamental chamber that contain $f$.  We denote the face corresponding to $f$ of the chamber labeled $w$ by $w(f)$.

\begin{thm}
$\Cm{W}$ can be constructed from $|W|$ copies of $\PG{W}$, labeled by the elements of $W$ and with the face $w(f)$ identified to $w'(f)$ whenever $w^{-1}w' \in  W_f$. The facet directly opposite $w$ through $w(f)$ is $ws_f$.
\label{t:glue}
\end{thm}

\begin{rem}
One may be tempted to think that whenever $w(f)$ is identified with $w'(f)$, the map between them is the restriction of the identity map between the chambers $w$ and $w'$. However, Davis et al.\ \cite[Section 8]{djs2} show that this is not the case, and compute the actual gluing maps between faces. For this reason they call the elements $s_f$ ``mock reflections.''  The gluing map may also be computed by applying the theorem above to subfaces of $f$.
\end{rem}

\begin{proof}
Since the $W$-action is type preserving, a chamber $w$ contains a face $f$ if and only if $w$ preserves $f$.  Recall that $W_f$ is generated by reflections in facets that contain $f$. Thus $f$ is contained only in chambers whose elements correspond to $W_f$. The chamber that lies directly across $f$ from the fundamental chamber corresponds to the longest word in $W_f$. Minimal blow ups of $\C{W}$ resolve nonnormal crossings, so $W_f$ is isomorphic to $(\Z/2\Z)^d$, where $F$ has codimension $d$. Thus the longest word in $W_f$ is the product of  generators $s_{f_i}$.

For every subspace $\al \in \Min{W}$ and every $w\in W$, the subspace $w(\al)$ is also in $\Min{W}$.  Thus we may extend the adjacency relation to chambers other than the fundamental chamber analogously.  Since the $W$-action preserves containment, a face $w'(f)$ is identified with $w(f)$ if and only if $w^{-1}w' \in W_f$.  Similarly, $w$ respects reflection across $F$ so the chamber directly across $w(f)$ from $w$ is $ws_f$.
\end{proof}

\begin{figure}[h]
\includegraphics{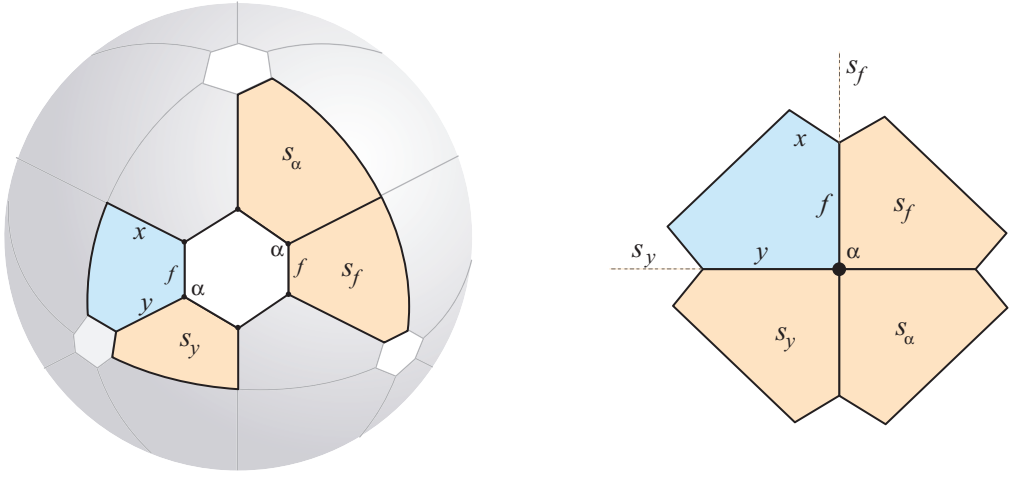}
\caption{Reflections locally around $\Cm{A_3}$.}
\label{f:gluing}
\end{figure}

\begin{exmp}
Consider the Coxeter group $A_3$.  Denote two facets of the fundamental chamber of $\C{A_3}$ by $x$ and $y$, whose reflections have the property that $(s_xs_y)^3=1$.  Note that $\Cm{A_3}$ is tiled by $24$ copies of the associahedron $\Pol (A_3)$. Let $f$ be the facet adjacent to $x$ and $y$ in the fundamental chamber of $\Cm{A_3}$ and let $\al$ be the intersection of $y$ and $f$, as in Figure~\ref{f:gluing}. Then $s_f=s_x s_y s_x$, and the fundamental chamber meets $s_y s_x$, $s_x s_y s_x$, and $s_y$ at $\al$.  If we travel directly across $\al$ from the fundamental chamber, we arrive in $s_{\al} = s_x s_y s_x \cdot s_y = s_y s_x$.
\end{exmp}

%%%%%%%%%%%%%%%%%%%%%%%%%%%%%%%%%%%%%%%%%%%%%%%%%%%%%%%%%%%%%%%%%%%%%%%%%%%%%%%%%%%%%
%
%                Geometry of Minimal Blow-ups
%
%%%%%%%%%%%%%%%%%%%%%%%%%%%%%%%%%%%%%%%%%%%%%%%%%%%%%%%%%%%%%%%%%%%%%%%%%%%%%%%%%%%%%

\newpage

\section{Geometry of Minimal Blow-ups}
\label{s:geometry}

\subsection{}
\label{ss:blowup}

One of our objectives is to describe the geometric structures of $\Min{W}$ before and after blow-ups.  This final section proves Theorem~\ref{t:blow_ups} which describes $\Cm{W}$ seen from the viewpoint of $\C W$.

We recall elementary notions of local structures, along with fixing notation:  The tangent space of a variety $V$ at $p$ is denoted $\TS{p}{V}$.   For a Coxeter complex $\C{W}$, the tangent space has a natural Euclidean geometry which it inherits from the embedding of $\C{W}$ in $\R^n$
(with the hyperbolic simplicial Coxeter groups being viewed as acting on the hyperboloid model inside $\R^n$).
%The tangent spaces of hyperbolic simplicial Coxeter groups also have this geometry since they can be viewed as acting on the hyperboloid model inside $\R^n$.
Two nonzero subspaces of $\TS{p}{V}$ are \emph{perpendicular} if each vector in the first is perpendicular to each vector in the second, under the Euclidean geometry of $\TS{p}{V}$.  The \emph{tangent bundle} of a variety $V$ on a subvariety $U$ is
$$\TS{U}{V} = \{(p,v)\suchthat p\in U, v\in \TS{p}{V} \}.$$
If $U = V$, we write $\TS{}{V}.$  The \emph{normal space} of $U$ at $p$ is
$$\NS{p}{U} = \{ v \suchthat v\in \TS{p}{V}, v \perp \TS{p}{U} \}$$
and the \emph{normal bundle} of $U$ at a subvariety $W\subset U$ is
$$\NS{W}{U} = \{ (p,v) \suchthat p\in U, v\in \NS{p}{V}\}.$$
If $W=U$, we write $\NS{}{U}$.

\begin{defn} \cite{har}
The \emph{blow-up} of a variety $V$ along a codimension $k$ intersection $\al$ of hyperplanes is the closure of $\{(x,f(x))\suchthat x\in V\}$ in $V \times \Pj^{k-1}$. The function $f:V\to \Pj^{k-1}$ is defined by $f:p \mapsto [f_1(p):f_2(p):\cdots :f_k(p)]$, where the $f_i$ define hyperplanes of $\Hs \al$ whose intersection is $\al$.
\end{defn}

We denote the blow-up of $V$ along $\al$ by $V_{\# \al}$.  There is a natural projection map
$$\pi: V_{\# \al} \to V : (x,y) \mapsto x$$
which is an isomorphism on $V-\al$.  The hyperplanes of $V_{\# \al}$ are the closures $\pi^{-1}(h-\al)$ for each hyperplane $h$ of $V$ and one additional hyperplane $\pi^{-1}(\al)$. Thus $V - \al$ and $V_{\# \al} - \pi^{-1}(\al)$ are isomorphic not only as varieties but as cellulations.\footnote{We give each hyperplane $h$ of $V_{\# \al}$ the same name as its projection $\pi(h)$.}  The hyperplane $\al$ of $V_{\# \al}$ has a natural identification with the projectified normal bundle of $\al$ in $V$. The intersection of a hyperplane $h$ with $\al$ is the part of $\al$ that corresponds to $\TS{\al}{h} \subset \NS{}{\al}$.

%%%%%%%%%%%%%%%%%%%%%%%%%%%%%%%%%%%%%%%%%%%%%%%%%%%%%%%%%%%%%%%%%%%%%%%%%%%%%%%%%%%%%

\subsection{}

A arrangement of hyperplanes of a variety $V$ cut $V$ into regions.  We say that the hyperplanes give a  \emph{cellulation} of $V$.  Two cellulations are equivalent if there is a hyperplane-preserving isomorphism between  the two varieties.  Let $\al$ be an intersection of hyperplanes. We say that hyperplanes $h_i$ \emph{cellulate}  $\al$ to mean the intersections $h_i \cap \al$ give a cellulation of $\al$, denoted by $\C{\al}$.  The notation $\C{\al}$ will always refer to the cellulation of $\al$ in the original complex, rather than its image in subsequent  blow-ups.  Let $\Min{\al}$ denote the minimal building set of $\C{\al}$, and let $\Cm{\al}$ denote the blow-up of  the minimal building set of $\al$.

\begin{thm}
\label{t:blow_ups}
Let $\C{W}$ be the Coxeter complex of a simplicial Coxeter group $W$ and let $\al \in \Min{W}$. The blow-up of \,$\al$\ in $\Cm{W}$ is equivalent to the product
$$ \Cm{\al} \times \Cpm{W_\al}.$$
\end{thm}

\begin{figure}[h]
\includegraphics{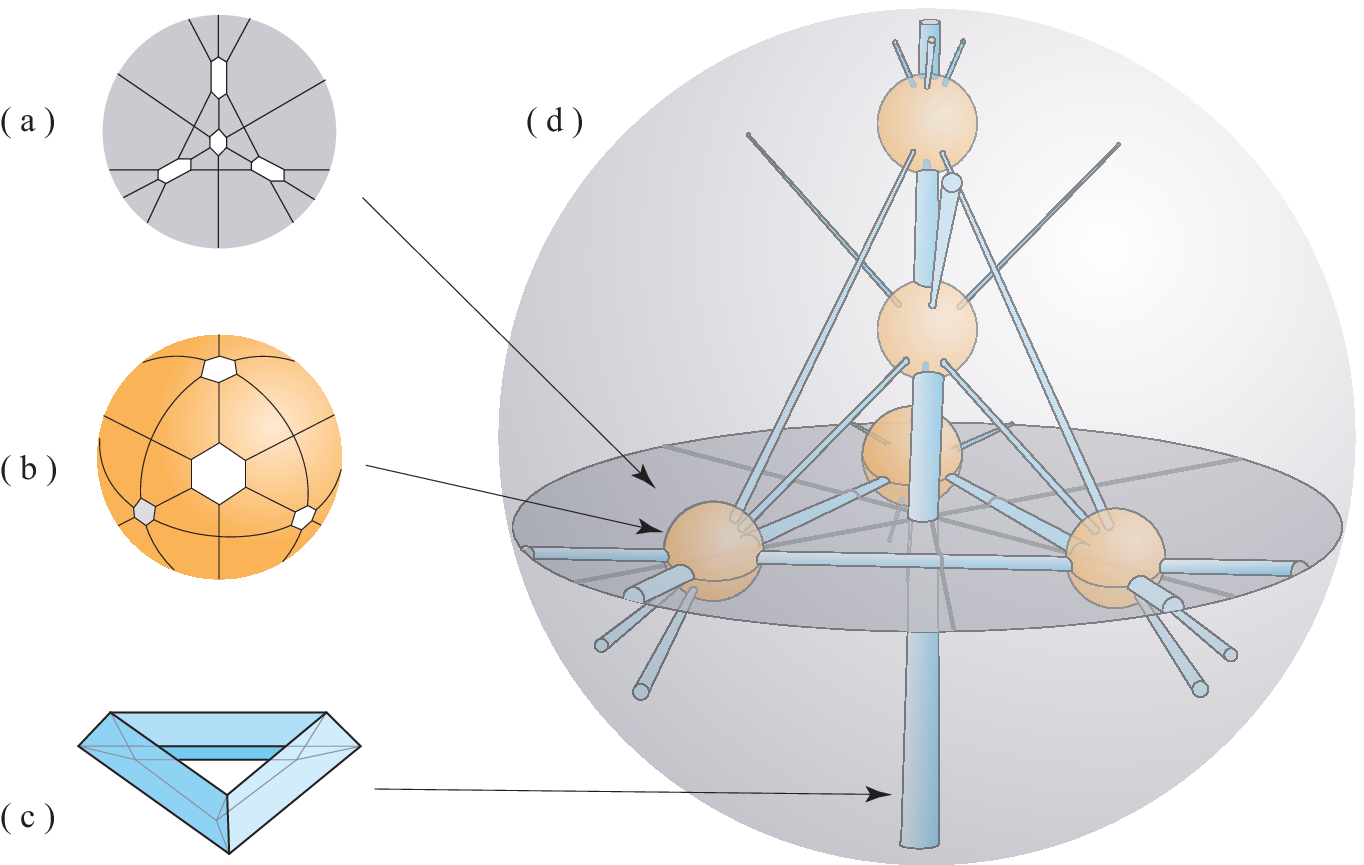}
\caption{The projective Coxeter complex (d) $\Cpm{A_4}$ along with components (a) $\Cpm{A_3}$, (b) $\Cm{A_3}$ and (c) $\Cpm{A_2} \times \Cpm{A_2}$.} \label{f:a4}
\end{figure}

\begin{exmp}
There are $2 \binom{n+1}{n-k}$ dimension $k$ elements of $\Min{A_n}$.  Each of these elements become $\Cm{A_{k+1}} \times \Cpm{A_{n-k-1}}$ in $\Cm{A_n}$.  Figure~\ref{f:a4}(d) shows the projective Coxeter complex $\Cpm{A_4}$ after minimal blow-ups. This is the Deligne-Knudsen-Mumford compactification \M{6} of the real moduli space of curves with six marked points.  It is the real projective sphere $\R\Pj^3$ with five points and ten lines blown-up.  Each of the five blown-up points are $\Cpm{A_3}$, shown in Figure~\ref{f:a4}(b) as $\Cm{A_4}$ before projecting through the antipodal map.  Each of the ten lines, each line defined by two distinct points in $\Min{A_4}$, becomes $\Cpm{A_2} \times \Cpm{A_2}$, a $2$-torus depicted in Figure~\ref{f:a4}(c).  Note that there are also ten codimension $1$ subspaces $\Cpm{A_3}$ pictured in Figure~\ref{f:a4}(a), defined by three distinct points in $\Min{A_4}$.
\label{e:a4}
\end{exmp}

\begin{rem}
Lemmas \ref{l:uniqreduce} and \ref{l:cellstab} are enough to provide the results of Theorem~\ref{t:blow_ups} for the maximal blow-up of $\C{W}$.
\end{rem}

\begin{rem}
Extensions of these results to configuration spaces are given in \cite[Section 3]{small}.
\end{rem}

%%%%%%%%%%%%%%%%%%%%%%%%%%%%%%%%%%%%%%%%%%%%%%%%%%%%%%%%%%%%%%%%%%%%%%%%%%%%%%%%%%%%%

\subsection{}

The proof of Theorem~\ref{t:blow_ups} requires two definitions and four preliminary lemmas.

\begin{defn}
Let $\be$ and $\ga$ be intersections of hyperplanes in a cellulation of $V$.  We say that $\be$ is \emph{strongly perpendicular} to $\ga$ and write $\be \stp \ga$ if for all $p$ in $\be \cap \ga$, all three of the following subspaces span $\TS{p}{V}$ \emph{and} any two of them are perpendicular: \newline
\indent (1) \ $\TS{p}{\be\cap \ga}$, \ (2) \ $\NS{p}{\be}$, \ and \ (3) $\ \NS{p}{\ga}$.
\end{defn}

Note that this directly implies that $\TS{p}{\be}$ is the span of $\TS{p}{\be\cap \ga}\cup \NS{p}{\ga}$.  For an intersection of hyperplanes $\be$, the normal space $\NS{p}{\be}$ is the span of the normal spaces of the elements of $\Hs \be$ at $p$; if $\be$ contains $\ga$, then $\NS{p}{\ga}$ contains $\NS{p}{\be}$.  This shows immediately that if $\Hs \be$ reduces to $\Hs \be_1 \sqcup \Hs \be_2$, then $\be_1 \stp \be_2$.

%%%%%%%%%%%%%%%%%%%%%%  LEMMA    %%%%%%%%%%%%%%%%%%

\begin{lem}
\label{l:uniqreduce} For every intersection of hyperplanes $\be$ in a Coxeter complex $\C{W}$, the set $\Hs \be$ has a unique maximal decomposition $\Hs \be=\Hs \be_1 \sqcup \cdots \sqcup \Hs \be_k$ where
\begin{enumerate}
\item each $\Hs \be_i$ is irreducible,
\item $\be_i \stp \be_j$ for all $i \neq j$, and
\item $\bigcap_{i\in S} \be_i$ properly contains $\be$ for any proper subset $S \subset \{ 1, 2, \ldots, k\}$.
\end{enumerate}
\end{lem}

\begin{proof}
If the normal spaces of two hyperplanes $h_1, h_2$ of $\Hs \be$ are not perpendicular, write $h_1 \sim h_2$. Then $\sim$ is a symmetric, reflexive relation on $\Hs \be$. Let $\approx$ be the unique smallest equivalence relation containing $\sim$ as a subset of $\Hs \be \times \Hs \be$.

No two hyperplanes $h_1 \sim h_2$ can be separated by any reduction of $\Hs \be$.  To prove this, suppose they could, and let $\Hs \be$ reduce to $\Hs \be_1 \sqcup \Hs \be_2$ with $h_1$ in $\be_1$ and $h_2$ in $\be_2$.  Then since $W_\be$ is a Coxeter group, the reflection of $h_1$ across $h_2$ must be in $\Hs \be$. By hypothesis, the resulting hyperplane must contain $\be_1$ or $\be_2$. The former implies that $\be_1 \subset h_2$ and the latter implies $\be_2 \subset h_1$, yielding a contradiction.  Since $\approx$ is the smallest transitive relation containing $\sim$, the hyperplanes $h_1, h_2$ cannot be separated whenever $h_1 \approx h_2$.

However, if $\approx$ partitions $\Hs \be$ into at least two classes, then we may separate $\Hs \be$ into $H_1 \sqcup H_2$ with each partition contained in either $H_1$ or $H_2$.  Clearly $\bigcap H_1 \cap \bigcap H_2= \be$.  To verify that $\Hs (\bigcap H_i)=H_i$, note that no element $h_1$ of $H_1$ may contain $\bigcap H_2$. If it does, then for all $p$ in $\be$, we have $\NS{p}{h_1}$ contained in $\NS{p}{\bigcap H_2}$, and thus in the span of $\{\NS{p}{h_2}\}$ for $h_2$ in $H_2$. This violates the pairwise perpendicularity in our choice of $H_1, H_2$.  Thus the equivalence relation $\approx$ partitions $\Hs \be$ into a unique maximal decomposition and therefore the $\Hs \be_i$'s are irreducible. Also, no proper subset $S$ of the $\be_i$ can intersect in exactly $\be$, since then $\cup \Hs \be_i$ for $\be_i$ in $S$ would be $\Hs \be$.  But $\Hs \be$ must reduce to $\cup \Hs \be_i$ for $\be_i \in S$ and $\cup \Hs \be_j$ for $\be_j \notin S$ by the argument above.

When $\Hs \ga$ reduces to $\Hs \ga_1 \sqcup \Hs \ga_2$, we have $\ga_1 \stp \ga_2$. Furthermore, since $\NS{p}{\ga_1}$ is the span of the normal spaces of $\Hs \ga_1$, and $\Hs \ga_1 \subset \Hs \ga$, then for any $\ga_3 \stp \ga$ (with nonempty intersection), it follows that $\ga_3 \stp \ga_1$. Thus by induction, $\be_i \stp \be_j$ for $i\neq j$.
\end{proof}

%%%%%%%%%%%%%%%%%%%%%%%%%%%%%%%%%%%%%%%%%%%%%%%%%%%%%%%%%%%%%%%%%%%%%%%%%%%%%%%%%%%%%

\subsection{}

The following two lemmas describe the effect of a blow-up on a cellulation.  The first lemma combines several facts that follow directly from the definitions of hyperplanes and blow-ups.  Note that as we perform blow-ups of $\C{W}$, the set of hyperplanes that contain a given $\be$ may change.  However, $\Hs \be$ is always assumed to refer to the set of hyperplanes that contain $\be$ in $\C{W}$.

%%%%%%%%%%%%%%%%%%%%%%  LEMMA    %%%%%%%%%%%%%%%%%%
% warning! hard coded projective cellulations ahead!

\begin{lem}
\label{l:cellstab}
Let $\be$ be a subvariety of \,$V$ with cellulation $C_1$. Suppose the tangent spaces of the hyperplanes $\Hs \be$ cellulate the normal bundle at each point $p$ with cellulation $C_2$.
\begin{enumerate}
\item The subvariety $\be$ of $V_{\# \be}$ is a product $C_1 \times \Pj C_2$.
\item The tangent space $\TS{p}{V_{\# \be}}$ for $p\in \be$ retains a local Euclidean structure.  Roughly speaking, $n-1$ of the coordinate vectors are in $\TS{p}{\be}$, and the other is parallel to the $1$-dimensional subspace of $\NS{\pi(p)}{\be}$ that corresponds to $p$.
\item For each hyperplane $h$ of $V$ that meets $\be$ at a subvariety $\ga \neq \be$, the hyperplane $h$ of $V_{\# \be}$ meets $\be$ at $\ga \times \Pj C_2$.
\item For each hyperplane $h$ of $V$ that properly contains $\be$, the hyperplane $h$ of $V_{\# \be}$ meets $\be$ at $C_1 \times h'$, where $h'$ is the image of $\TS{p}{h}$ in $\Pj C_2$. Also $\be \stp h$.
\end{enumerate}
\end{lem}

%%%%%%%%%%%%%%%%%%%%%%  LEMMA    %%%%%%%%%%%%%%%%%%

\begin{lem}
\label{l:perpblow}
Let $\be$ be a subvariety of $V$ with cellulation $C_1$.  If $\be \stp \ga$ in $V$, then $\C{\be}$ in $V_{\# \ga}$ is equivalent to $(C_1)_{\#(\be\cap \ga)}$.
\end{lem}

\begin{proof}
The normal bundle  $\NS{\be\cap \ga}{\ga}$ is contained in $\TS{\be\cap \ga}{\be}$ since $\be \stp \ga$.  Thus $\NS{}{\be \cap \ga}$ and $\NS{}{\ga}$ have the same intersection with $\TS{}{\be}$. Since blow-ups replace a variety with its projectified normal bundle, the blow-ups along $\ga$ and $\be\cap \ga$ produce equivalent cellulations of $\be$.
\end{proof}

Finally we establish the tools that will allow us to change the order in which we blow up elements of $\Min{W}$.  The following definition and lemma give a class of orderings that produce the same cellulation as minimal blow-ups.

\begin{defn}
Given a variety $V$ with intersections of hyperplanes $\be, \ga$, the blow-ups along $\be$ and $\ga$ \emph{commute} if the cellulations $(V_{\# \be})_{\# \ga}$ and $(V_{\# \ga})_{\# \be}$ are equivalent, and the induced map on the hyperplanes preserves their labels.
\end{defn}

%%%%%%%%%%%%%%%%%%%%%%  LEMMA    %%%%%%%%%%%%%%%%%%

\begin{lem}
\label{l:reorder}
Let $x_1, x_2, \ldots, x_k$ be an ordering of the elements of $\Min{W}$ such that $i \leq j$ whenever $x_i$ is contained in $x_j$. Then blowing up $\C{W}$ along the $x_i$ in order gives a cellulation equivalent to $\Cm{W}$.  The induced map on the hyperplanes also preserves labels.
\end{lem}

\begin{proof}
First we verify that if $\be \stp \ga$, then the blow-ups along $\be$ and $\ga$ commute.  Since $\be \stp \ga$, the bundle $\NS{\be\cap \ga}{\be}$ is contained in $\TS{\be\cap \ga}{\ga}$ and $\NS{\be\cap \ga}{\ga}$ is contained in $\TS{\be\cap \ga}{\be}$. Define the maps $\pi_\be:V_{\# \be} \to V$ and $\pi_{\be\ga}:(V_{\# \be})_{\# \ga} \to V_{\# \be}$. Then $\pi_\be^{-1}(V-\ga)=\pi_\be^{-1}(V)-\pi_\be^{-1}(\ga)$ since $\NS{\be\cap \ga}{\be}\subset \TS{\be\cap \ga}{\ga}$.  Thus $(V_{\# \be})_{\# \ga}$ is the closure of $\pi_{\be\ga}(\pi_\be^{-1}(V-\ga))$, which is the closure of $\pi_{\be\ga}^{-1}(\pi_\be^{-1}(V-\ga-\be))$. Similar reasoning shows that $(V_{\# \ga})_{\# \be}$ is the closure of $\pi_{\ga\be}^{-1}(\pi_\ga^{-1}(V-\be-\ga))$. Since the $\pi$'s are isomorphisms on $V - \be -\ga$, we have a natural isomorphism between $(V_{\# \be})_{\# \ga}$ and $(V_{\# \ga})_{\# \be}$.

Now take $\be, \ga$ to be elements of $\Min{W}$ such that neither contains the other. By Lemma~\ref{l:minirred}, the arrangements $\Hs \be$ and $\Hs \ga$ are irreducible. Applying Lemma~\ref{l:uniqreduce} shows that if $\Hs (\be \cap \ga)$ is reducible, then  $\be \stp \ga$ and $\Hs \be \sqcup \Hs \ga$ is the unique reduction.

If $\Hs (\be \cap \ga)$ is irreducible, then $\be\cap \ga$ is in $\Min{W}$. After the blow-up along $\be \cap \ga$, the resulting spaces $\be$ and $\ga$ do not intersect by Lemma~\ref{l:cellstab}, and thus (vacuously) $\be \stp \ga$.  If $\be\cap \ga$ is not in $\Min{W}$, then $\be \stp \ga$. In either case, the blow-ups along $\be$ and $\ga$ commute.  Thus we may transpose any two elements that do not contain each other in the ordering of $\Min{W}$ and get an equivalent cellulation (with matching hyperplane labels) after blowing up all of $\Min{W}$. Repeating this procedure proves the statement of the lemma.
\end{proof}

%%%%%%%%%%%%%%%%%%%%  BIG PROOF   %%%%%%%%%%%%%%%%%%%%%%%%%%%%%%%%%%%%%%%%%%%%%%%%%%%%%%%%%%%

\subsection{}

We have now assembled all the lemmas needed for the proof of the theorem.

\begin{proof}[Proof of Theorem~\ref{t:blow_ups}]

We begin by applying Lemma~\ref{l:reorder}.  Divide the elements of $\Min{W}$ into three sets:
\begin{enumerate}
\item $\{\al\}$,
\item $X = \{\be : \be \not\supset \al\}$, and
\item $Y = \{\be: \be \subset \al\}$.
\end{enumerate}

We reorder the elements of $\Min{W}$ as follows:  First we blow up the elements of $X$, ordered by the dimension of $\be\cap \al$, followed by blowing up along $\al$.  Finally blow up the elements of $Y$ in order of dimension, as usual.  Note that this is a valid application of Lemma~\ref{l:reorder}, since if $\be$ contains $\ga$, then $\be\cap \al$ contains $\ga \cap \al$.

We next produce a bijection $\phi$ between the set $X'$ of elements $x_i$ in $X$ that intersect $\al$ in $(\cdots((\C{W})_{\#x_1})_{\#x_2}\cdots)_{\#x_{i-1}}$ and the elements of $\Min{\al}$ in $\C{W}$.  We show that the map $\phi: X' \to \Min{\al}: \be \mapsto \be \cap \al$ is a bijection, and that blowing up the elements of $X$ has the same effect on the cellulation of $\al$ as blowing up the elements of $\Min{\al}$.

\begin{enumerate}
\item Suppose $\be \in X'$ and $\be \subset \al$, and thus $\Hs \al \subset \Hs \be$. Since $\be$ is in $\Min{W}$, the group $W_\be$ is an irreducible spherical Coxeter group by Lemma~\ref{l:irreducible}, and the arrangement $\Hs \al$ is irreducible in $\Hs \be$ by Lemma~\ref{l:minirred}.
% One may directly verify for all spherical Coxeter groups $W_\be$ that the elements of $\Hs \be$ intersect $\al$ in an irreducible arrangement.
For all spherical Coxeter groups $W_\be$, the elements of $\Hs \be$ intersect $\al$ in an irreducible arrangement.\footnote{This can be checked by hand for the simpler cases, and a detailed decomposition of types $A_n$, $B_n$ and $D_n$ is given in \cite[Section 5]{small}.  For the larger complexes ($E_6$, $E_7$, $E_8$), we can exploit the appearance of $A_n$ as a subgroup.}
Therefore $\be=\be \cap \al \in \Min{\al}$.

\item Suppose $\be \in X'$ and $\be \not\subset \al$, then $\be \cap \al$ is not in $\Min{W}$, thus $\Hs (\be \cap \al)$ reduces. Lemma~\ref{l:uniqreduce} guarantees that $\be \stp \al$. Thus by Lemma~\ref{l:perpblow}, blowing up $\be$ is equivalent to blowing up $\be\cap \al$ in the cellulation of $\al$.

\item We now produce an function $\psi: \Min{\al} \to X'$ that will be the inverse to $\phi$.  For $\be \in \Min{\al}$, either $\be \in \Min{W}$ or $\Hs \be$ is reducible. If $\be \in \Min{W}$, then let $\psi(\be) = \be$.  If not, then $\Hs \be$ must reduce to $\Hs \al_0 \sqcup \Hs \al_1 \sqcup \cdots \sqcup \Hs \al_m$.  Without loss of generality, assume $\al$ contains $\al_0$.

Since $\NS{\be}{\al_i}$ is contained in $\TS{\be}{\al_0}$ for $i\neq 0$, the normal spaces of the elements of $\Hs \be-\Hs \al_0$ are the same in $\al$ as they are in $\C{W}$. Thus in $\al$, we know that $\al_i\stp \al_j$ for $i, j \neq 0, i\neq j$. Furthermore, the normal space $\NS{p}{\al_0}$ in $\al$ is a subset of $\NS{p}{\al_0}$ in $V$.  Thus if $\NS{p}{\al_0}$ in $\al$ is nonzero, it is perpendicular to each $\NS{p}{\al_i}$ in $\al$. Thus the set of hyperplanes of $\al$ induced by $\Hs \be-\Hs \al$ reduces to the disjoint union induced by $(\Hs \al_0 -\Hs \al )\sqcup \Hs \al_1 \sqcup \cdots \sqcup \Hs \al_m$. To satisfy the hypothesis that $\be \in \Min{\al}$, it is necessary that $\al_0=\al$ and $m=1$.  Thus we define $\psi(\be)=\al_1$.   It is straightforward to check that $\psi$ is the inverse of $\phi$, so $\phi$ is a bijection.

\item By our choice of ordering, the elements of $X'$ are blown up in the same order as elements of $\Min{\al}$ under minimal blow-ups. Furthermore, the subvariety $\al$ has equivalent cellulations in the blow-up along $\phi(\be)$ and in the blow-up along $\be$. This follows trivially if $\be\in \al$, and from Lemma~\ref{l:perpblow} if not.

\end{enumerate}

Thus, after blowing up all the elements of $X$ in $\C{W}$, the cellulation of $\al$ is equivalent to $\Cm{\al}$.  By Lemma~\ref{l:cellstab}, the result after blowing up $\al$ is equivalent to $\Cm{\al} \times \Cp{W_\al}$. Furthermore, for each element $y\in Y$, we have $y \stp \al$ and $y \cap \al=\Cm{\al} \times y'$, where $y'$ is the image of $y$ in $\Cp{W_\al}$. Since the elements of $Y$ are ordered by dimension, they are also ordered by their dimension in $\C{W_\al}$. Lemma~\ref{l:perpblow} guarantees that blowing up the elements of $Y$ produces a cellulation of $\al$ equivalent to $\Cm{\al} \times \Cpm{W_\al}$.
\end{proof}

%%%end of proof

\begin{ack}
We thank Mike Davis, Rick Scott and Jim Stasheff for advice and clarifications.  We also thank Zan Armstrong, Eric Engler, Ananda Leininger and Michael Manapat for numerous discussions.
\end{ack}

%%%%%%%%%%%%%%%%%%%%%%%%%%%%%%%%%%%%%%%%%%%%%%%%%%%%%%%%%%%%%%%%%%%%%%%%%%%%%%%%%%%%%
%
%                  REFERENCES
%
%%%%%%%%%%%%%%%%%%%%%%%%%%%%%%%%%%%%%%%%%%%%%%%%%%%%%%%%%%%%%%%%%%%%%%%%%%%%%%%%%%%%%
\bibliographystyle{amsplain}

\begin{thebibliography}{ADM9}
\baselineskip=12pt

\bibitem[1]{small} S.\ Armstrong, M.\ Carr, S.\ Devadoss, E.\ Engler, A.\ Leininger, and M.\ Manapat, Point configurations and Coxeter operads, preprint 2004.

\bibitem[2]{bt} R.\ Bott and C.\ Taubes, On the self-linking of knots, \emph{J.\ Math.\ Phys.} {\bf 35} (1994), 5247-5287.

\bibitem[3]{bou} N.\ Bourbaki, \emph{Lie Groups and Lie Algebras: Chapters 4-6,} Springer-Verlag, Berlin, 2002.

\bibitem[4]{bro} K.\ S.\ Brown, \emph{Buildings,} Springer-Verlag, New York, 1989.

\bibitem[5]{djs} M.\ Davis, T.\ Januszkiewicz, R.\ Scott, Nonpositive curvature of blowups, \emph{Selecta Math.} {\bf 4} (1998), 491 - 547.

\bibitem[6]{djs2} M.\ Davis, T.\ Januszkiewicz, R.\ Scott, Fundamental groups of minimal blow-ups, \emph{Adv. Math.} {\bf 177} (2003), 115-179.

\bibitem[7]{dp} C.\ De Concini and C.\ Procesi, Wonderful models of subspace arrangements, \emph{Selecta Math.} {\bf 1} (1995), 459-494.

\bibitem[8]{dev} S.\ Devadoss, Tessellations of moduli spaces and the mosaic operad,  {\em Contemp.\ Math.} {\bf 239} (1999), 91-114.

\bibitem[9]{dev2} S.\ Devadoss, A space of cyclohedra, \emph{Disc.\ Comp.\ Geom.} {\bf 29} (2003), 61-75.

\bibitem[10]{dev3} S.\ Devadoss, Combinatorial equivalence of real moduli spaces, \emph{Notices Amer.\ Math.\ Soc.} {\bf 51} (2004), 620-628.

\bibitem[11]{fz} S.\ Fomin, A.\ Zelevinsky, $Y$-systems and generalized associahedra, \emph{Ann. Math.} {\bf 158} (2003), 977-1018.

\bibitem[12]{har} R.\ Hartshorne, \emph{Algebraic Geometry,} Springer-Verlag, New York, 1977.

\bibitem[13]{kap1} M.\ M.\ Kapranov, The permutoassociahedron, MacLane's coherence theorem, and asymptotic zones for the $KZ$ equation, \emph{J.\ Pure Appl.\ Alg.} {\bf 85} (1993), 119-142.

\bibitem[14]{lee} C.\ Lee, The associahedron and triangulations of the $n$-gon, \emph{Euro.\ J.\ Combin.} {\bf 10} (1989), 551-560.

\bibitem[15]{sta} J.\ D.\ Stasheff, Homotopy associativity of $H$-spaces, \emph{Trans.\ Amer.\ Math.\ Soc.} {\bf 108} (1963), 275-292.

\bibitem[16]{sta2} J.\ D.\ Stasheff (Appendix B coauthored with S.\ Shnider), From operads to ``physically'' inspired theories, \emph{Contemp.\ Math.} {\bf 202} (1997), 53-81.

\end{thebibliography}

\end{document}